\documentclass[centertags]{article}
\input epsf
\usepackage{amsmath}
\usepackage{amsthm}
\usepackage{amscd}
\usepackage{amssymb}
\usepackage{verbatim}

\title{Rotation Numbers and Instability Sets}
\author{John Franks\thanks{Supported in part by NSF grant DMS0099640.
This article is the written version of an invited address at the
January, 2002 AMS meeting in San Diego, California.}}

\newtheorem{thm}{Theorem}[section] 
\newtheorem{conj}[thm]{Conjecture}
\newtheorem{quest}[thm]{Question}
\newtheorem{cor}[thm]{Corollary}
 
\newtheorem{prop}[thm]{Proposition}
\newtheorem{defn}[thm]{Definition}

\newcommand{\R}{{\mathbb R}}
\newcommand{\Z}{{\mathbb Z}}
\newcommand{\A}{{\mathbb A}}
\newcommand{\tiA}{\tilde {\mathbb A}}
\newcommand{\T}{{\mathbb T}}
\newcommand{\Q}{{\mathbb Q}}
\newcommand{\I}{I}

\DeclareMathOperator{\dist}{dist}

\DeclareMathOperator{\Per}{Per}

\DeclareMathOperator{\Diff}{Diff}

\begin{document}
\maketitle
\begin{abstract}
Translation and rotation numbers have played an interesting and
important role in the qualitative description of various dynamical
systems.  In this exposition we are especially interested in
applications which lead to proofs of periodic motions in various kinds
of dynamics on the annulus.  The applications include billiards and
geodesic flows.

Going beyond this simple qualitative invariant in the study of 
the dynamics of area preserving annulus maps, G.D. Birkhoff was led to
the concept of ``regions of instability'' for twist maps.
We discuss the closely related notion of instability sets for a
generic area preserving surface diffeomorphism and develop their properties.
\end{abstract}

\section{Introduction}
The concepts of translation number and rotation number were introduced
by Henri Poincar\'e and have had a rich history. We describe several
of many manifestations of these concepts as they apply to the
qualitative description of various dynamical systems. We are
especially interested in applications which lead to proofs of periodic
motions in various kinds of area preserving dynamics on the annulus
and $S^2$.  Applications include billiards and geodesic flows.  This
is the content of sections \S2 to \S5.

In pursuit of the qualitative analysis of such dynamics, G.D. Birkhoff
was led to the concept of ``regions of instability'' for twist maps on
an annulus.  We discuss the closely related notion of instability
sets for a generic area preserving surface diffeomorphism and describe
some of their properties.

Guiding this is the deep theory of twist
maps of the annulus where we have a large body of knowledge, thanks to
the pioneering work of G.D. Birkhoff and subsequent work by Arnold,
Kolmagorov, Herman, Mather, and Moser.  However, the definition of an
annular twist map depends on a distinguished choice of co-ordinates on
the annulus. In a local setting, e.g. in the neighborhood of a generic
fixed point, one can always obtain these co-ordinates, but they are
generally unavailable in the larger global context. It would be nice
to have the same or analogous results in this global context.

In sections \S~\ref{sec: instab}, \S~\ref{sec: properties} and 
\S~\ref{sec: connections} we describe a residual class of $C^r$ 
area preserving surface diffeomorphisms for $r \ge 16$ and describe
some of their qualitative properties. 

In \S9 we discuss some of the intriguing open questions related to
instability sets.

\section{Circle Homeomorphisms}
The simplest instance of rotation number is its use as an invariant
for orientation preserving circle homeomorphisms.  Define
the circle $\T^1$ to be the quotient, $\R/\Z.$
As a quotient $\T^1$ inherits both a group structure and a topology.
There is a natural projection $\pi: \R \to \T^1$ defined by
$\pi(x) = x + \Z$ which is both a (universal) covering map and
a group homomorphism.

\begin{defn}
Let $f : \T^1 \to \T^1$ be an orientation preserving homeomorphism 
of the circle.  A homeomorphism $F: \R \to \R$ is a {\em lift} of $f$ provided
\[
f( \pi(t)) = \pi( F(t)).
\]
\end{defn}

If $F$  is a lift of $f$ it is not difficult to see that
$G: \R \to \R$ is also a lift of $f$ if and only if $F(t) = G(t) + n$ for
some $n \in \Z.$

\begin{defn}
Let $f : \T^1 \to \T^1$ be an orientation preserving homeomorphism 
and let $F: \R \to \R$ be a lift.  The {\em translation number\footnote{It 
is common to see this quantity referred to as the 
{\it rotation number}, however, the term translation number,
while less commonly used, is more accurately descriptive.
More importantly, there is a
need to distinguish the translation number from its projection in $\T^1,$ 
which is properly called the rotation number (see 
Definition~(\ref{defn: rotation number})).}}
of a point $x \in \R$ under $F$ is
\[
\tau(x,F) = \lim_{n \to \infty} \frac{ F^n(x) - x}{n}.
\]
\end{defn}

Notice that the value of $\tau(x, F)$ is the average amount the point $x$ is
translated under iteration by $F$.
We summarize the elementary properties of the translation number of
a lift $F$ in the following proposition (See section 11.1 of \cite{KH}
for a proof).

\begin{prop}\label{prop: lift properties}
Let $f : \T^1 \to \T^1$ be an orientation preserving homeomorphism 
and let $F: \R \to \R$ be a lift.
\begin{itemize} 
\item
The limit $\tau(x, F)$ exists for all $x \in \R.$
\item
$\tau(x, F^n) = n \tau(x, F)$ for all $n \in \Z.$
\item
$\tau(x, F + m) = \tau(x, F) + m$ for $m \in \Z.$
\item
$\tau(F) := \tau(x, F)$  is independent of $x$.
\end{itemize}
\end{prop}

More striking is the fact that if $\tau(F)$ is rational, then there
are of periodic orbits.

\begin{thm}\label{thm: rat => periodic}
Let $f : \T^1 \to \T^1$ be an orientation preserving homeomorphism 
and let $F: \R \to \R$ be a lift.
\begin{itemize} 
\item
$\tau(F) = 0$  if and only if $F$ has a fixed point.
\item
$f$ has a fixed point if and only if $\tau(F) \in \Z$.
\item
$f$ has a periodic point of period $q$ if and only if
$\tau(F) = p/q \in \Q$, where $p/q$ is in lowest terms.
\end{itemize}
\end{thm}

If $F$ and $G$ are two different lifts of $f: \T^1 \to \T^1$, then
they differ by a constant integer, i.e. there is an $m \in \Z$
such that $F(t) = G(t) + m$ for all $t \in \R$.  From this and the
proposition above it follows that $\tau(x, F) = \tau(x, G) + m$.
As a consequence if we reduce $\tau(x, F)$ modulo 1, or, more 
precisely, consider its projection in $\T^1 := \R/\Z$ the result
is independent of the choice of the lift $F$ and depends only
on the original homeomorphisms $f.$  

\begin{defn}\label{defn: rotation number}
Let $f : \T^1 \to \T^1$ be an orientation preserving homeomorphism 
and let $F: \R \to \R$ be a lift.  The {\em rotation number}
of $f$ is the element of $\T^1$ defined by
\[
\rho(f) = \pi( \tau(x, F))
\]
where $\pi : \R \to \T^1$ is the natural projection, $F$ is any 
lift of $f$ and $x$ is any point of $\R.$
\end{defn}

From Proposition~\ref{prop: lift properties} we immediately obtain the
following

\begin{prop}\label{prop: T^1 properties}
Let $f : \T^1 \to \T^1$ be an orientation preserving homeomorphism.
\begin{itemize} 
\item 
$\rho(f^n) = n \rho(F)$ for all $n \in \Z.$
\item
$f$ has a periodic point of period $q$ if and only if
        $\rho(f) = (p/q + \Z) \in \Q/\Z$ where $p/q$ is in lowest terms.
\end{itemize}
\end{prop}

For the many important dynamical systems which preserve a measure
there is often an important relation between the measure
and the rotation numbers. In the case of circle homomorphisms this 
relation has an especially simple form which we will now describe.  

\begin{defn}
Let $f : \T^1 \to \T^1$ be an orientation preserving homeomorphism 
and let $F: \R \to \R$ be a lift.  We define the 
{\em displacement function}  $\phi: \T^1 \to \R$ of $F$ by 
$\phi(x + \Z) = F(x) - x$.
\end{defn}

Note that if $x'$ is another choice for $x$, (i.e.
$x +\Z = x' +\Z$), then $x'  = x +m$
for some integer $m$ so $F(x) - x = F(x') - x'$ and 
$\phi$ is well defined.  For a point $(x +\Z) \in \T^1$ the value of 
$\phi$ can be thought of as amount this point
is displaced around the circle.  Having chosen a lift $F$ we have determined
this as a well defined element of $\R$, not just an element of $\R$
modulo $1.$

\begin{thm}\label{thm: T^1 mean}
Let $f : \T^1 \to \T^1$ be an orientation preserving homeomorphism 
which preserves a probability measure $\mu,$ and let $F$ be a lift
of $f$  and suppose $\phi: \T^1 \to \R$ is the displacement function
of $F$.  Then
\[
\tau(F) =  \int_{\T^1} \phi \ d\mu.
\]
\end{thm}

Notice that the left hand side of this equality is independent of $\mu$ 
and hence the right hand side is the same for all choices of $\mu$!

\begin{proof}
Here we use $z$ to denote a point of $\T^1$ instead
of the more cumbersome $x +\Z.$ 
According to the Birkhoff Ergodic Theorem (see \cite{KH}),
the function $\hat \phi$ defined by
\[
\hat \phi(z) = \lim_{n \to \infty} \frac{1}{n} \sum_{i=0}^{n-1} \phi( f^i(z))
\]
is integrable and satisfies
\[
\int_{\T^1} \phi \ d\mu = \int_{\T^1} \hat \phi \ d\mu.
\]
But if $z = \pi(x)$, then $f^i(z) = \pi(F^i(z))$ and
$\phi( f^i(z)) = F(F^{i}(x)) - F^{i}(x)  = F^{i+1}(x) - F^{i}(x).$
It follows that the telescoping sum
\[
\sum_{i=0}^{n-1} \phi( f^i(z)) 
= \sum_{i=0}^{n-1} \big (F^{i+1}(x) - F^{i}(x)\big )
= F^n(x) - x.
\]
Hence
\[
\hat \phi(z) = \lim_{n \to \infty} \frac{F^n(x) - x}{n} = \tau(F).
\]
So $\hat \phi$ is the constant function with value $\tau(F)$
and the result follows.
\end{proof}

\section{Annulus Homeomorphisms}

The annulus $\A = \T^1 \times \I$, where
$\I = [0,1],$
provides much richer ground for application of the
concept of rotation number.

Before exploring this, however, we motivate our interest in
the annulus by describing one of the many mechanical systems
which give rise to annulus maps.  Poincar\'e's original interest
in rotation numbers was related to the restricted three body
problem, but, a system much easier to describe is provided by
billiards on a strictly convex table with smooth boundary.

The ball is idealized as a point which moves on the table 
with uniform motion and
bounces off the edge in such a way that the angle of
incidence equals the angle of reflection.  The annulus associated with
this dynamical system is $\T^1 \times [0,\pi]$ where the circle $\T^1$
parametrizes the edge of the table and $\theta \in [0,\pi]$ is the
angle between the tangent to the table edge and the direction at which
the the ball leaves the edge at a bounce.  If $x$ is a point
of contact and $\theta \in (0,\pi)$ is an
angle of reflection, then the {\em billiard map} $F :
\A \to \A$ is defined by $F(x, \theta) = (x', \theta'),$ where $x'$ is
the next point of contact with the boundary and $\theta'$ is the next
angle of reflection.  Notice that by defining $F(x, 0) = x$
and $F(x, \pi) = x$ we can extend $F$ to a homeomorphism of all of $\A.$

To understand the billiard dynamics qualitatively (for example
the existence of periodic orbits) it suffices to understand 
this billiard map.  A crucial aid in this analysis comes from the
fact the map preserves an area.

\begin{thm}[Birkhoff]
If we define an area on $\A$ using the element of 
area $\sin \theta\ dx\ d\theta$, then this area is
preserved by the billiard map  $F : \A \to \A.$
\end{thm}

As in the case of a circle we can define the notions of
translation number and rotation number.

\begin{defn}
Let $\tiA = \R \times \I.$
A homeomorphism $F: \tiA \to \tiA$ is a {\em lift} of $f$ provided
\[
f(\bar \pi(x, y)) = \bar \pi( F(x,y)),
\]
where $\bar \pi: \tiA \to \A$ is given by $\bar \pi(x,y) = (\pi(x),y).$
\end{defn}

If $F$  is a lift of $f$ it is not difficult to see that
$G: \tiA \to \tiA$ is also a lift of $f$ if and only if
\[
F(x,y) = G(x,y) + (m, 0) 
\]
for some $m \in \Z.$

\begin{defn}\label{defn: A translation number}
Let $f : \A \to \A$ be an orientation preserving homeomorphism 
which preserves boundary components and let \\
$F: \tiA \to \tiA$ be a lift.  The {\em translation number}
of a point $w = (x,y) \in \tiA$ under $F$ is 
\[
\tau(w,F) = \lim_{n \to \infty} \frac{ p_1(F^n(w) - w)}{n},
\]
where $p_1: \tiA \to \R$ is projection on the first co-ordinate.
\end{defn}

We would like an analogue of Proposition~\ref{prop: lift properties},
i.e. we would like for the properties of the translation number
for circle homeomorphisms to hold for annulus homeomorphisms.
Unfortunately, much less is true for the annulus.  In particular
$\tau(w, F)$ may fail to exist for many $w \in \tiA.$  Moreover,
when it does exist, it is only rarely independent of $w$.

\begin{prop}\label{prop: A lift properties}
Let $f : \A \to \A$ be an orientation preserving homeomorphism 
which preserves boundary components and let $F: \tiA \to \tiA$ be a lift.
\begin{itemize} 
\item
If $\tau(w, F)$ exists, then $\tau(w, F^n) = n \tau(w, F)$ for all $n \in \Z.$
\item
If $\tau(w, F)$ exists, then $\tau(w, F + (m,0)) = \tau(w, F) + m$ for $m \in \Z.$
\item
If $f$ preserves a finite Borel measure $\mu$ and $\tilde \mu$ is the
lift of this measure to $\tiA$, then $\tau(w, F)$ exists for all $w$
except a set of $\tilde \mu$ measure zero.
\end{itemize} 
\end{prop}

The proofs of the first two parts of this proposition are straightforward.
The final part is a consequence of the Birkhoff ergodic theorem and we
give a proof below in Theorem~\ref{thm: A mean translation}.

\begin{defn}\label{defn: A rotation number}
Let $f : \A \to \A$ be an orientation preserving homeomorphism 
and let $F: \tiA \to \tiA$ be a lift.  The {\em rotation number}
of $f$ is the element of $\T^1$ defined by
\[
\rho(f) = \pi( \tau(x, F))
\]
where $\pi : \R \to \T^1$ is the natural projection, $F$ is any 
lift of $f$ and $x$ is any point of $\R.$
\end{defn}

From Proposition~\ref{prop: A lift properties} we immediately obtain the
following:

\begin{prop}\label{prop: A rot properties}
Let $f : \A \to \A$ be an orientation preserving homeomorphism 
which preserves boundary components and let $F: \tiA \to \tiA$ be a lift.
\begin{itemize} 
\item
If $\tau(w, F)$ exists, then $\rho(\bar \pi(w), f)$ exists and is well
defined, i.e. independent of the lift $F$.
\item
If $\rho(x, f)$ exists $\rho(w, f^n) = n \rho(x, f)$ for all $n \in \Z.$
\item
If $f$ preserves a finite Borel measure $\mu$
then $\rho(x, f)$ exists for all $x \in \A$ except
a set of $\mu$ measure zero.
\end{itemize} 
\end{prop}

It is not hard to see that periodic orbits have rational translation
numbers.  But the issue of using translation numbers
or rotation numbers to prove the existence of periodic points is much
more subtle than in the case of the circle and that is the subject of
\S\ref{sec: P-B} below.  

The last item in Proposition~\ref{prop: A rot properties}, 
makes it very important to have a measure
preserving annulus homeomorphism and it is in this setting that we are
able to provide the most interesting results.

As in the case of $\T^1$ it is useful to define the displacement function
for a lift $F$.

\begin{defn}
Let $f : \A \to \A$ be an orientation preserving homeomorphism,
preserving boundary components, 
and let $F: \tiA \to \tiA$ be a lift.  We define the 
{\em displacement function}  $\phi: \A \to \R$ of $F$ by 
$\phi(\bar \pi(w)) = p_1(F(w) - w)$ for $w \in \tiA,$ 
where $p_1: \tiA \to \R$ is projection on the first co-ordinate.
\end{defn}

Note that if $w'$ satisfies $\bar \pi(w) = \bar \pi(w')$
then $w'  = w +(m,0)$
for some integer $m$ so $F(w) - w = F(w') - w'$ and 
$\phi$ is well defined.  
 
\begin{thm}\label{thm: A mean translation}
Let $f : \A \to \A$ be an orientation preserving homeomorphism,
preserving boundary components, 
and preserving a probability measure $\mu,$ and let $F$ be a lift
of $f$.  Then the rotation number $\rho(x, F)$ exists for
$\mu$-almost all $x$ and is a $\mu$ integrable
function. If $\phi: \A \to \R$ is the displacement function
of $F$, 
\[
\int_{\T^1} \rho(x,F) \ d\mu = \int_{\T^1} \phi \ d\mu.
\]
The value of these integrals will be denoted
by $\tau_\mu(F)$ and is called the {\em mean translation
number} of $F$.
\end{thm}

\begin{proof}
Again we appeal to the Birkhoff Ergodic Theorem (see \cite{KH}).
As in the proof of Theorem~\ref{thm: T^1 mean}, if $x = \bar \pi(w)$, then 
the telescoping sum
\[
\sum_{i=0}^{n-1} \phi( f^i(x)) 
= \sum_{i=0}^{n-1} p_1(F^{i+1}(w) - F^{i}(w))
= p_1(F^n(w) - w).
\]
So
\[
\tau(w,F) = \lim_{n \to \infty} \frac{p_1(F^n(w) - w)}{n}
= \lim_{n \to \infty} \frac{1}{n} \sum_{i=0}^{n-1} \phi( f^i(x)).
\]
Hence, according to the ergodic theorem, $\tau(w,F)$ exists
for $\mu$-almost all $x$ and is a $\mu$  integrable.  Moreover
it satisfies
\[
\tau_\mu(F) = \int_{\T^1} \tau \ d\mu = \int_{\T^1}  \phi \ d\mu.
\]
\end{proof}

\section{The Poincar\'e-Birkhoff Theorem}\label{sec: P-B}

In this section we describe several theorems analogous to 
Theorem~\ref{thm: rat => periodic}, in the sense that from a hypothesis on a
translation number one is able to conclude the existence of
a periodic orbit.  The most famous of these is the 
Poincar\'e-Birkhoff Theorem.  This was conjectured by Poincar\'e
and he was able to prove some special cases.  His motivation
was to prove the existence of some periodic motions of the restricted
three-body problem (this is the Newtonian motion of three point
masses, one of which has zero mass).

G.~D.~Birkhoff who succeeded in proving the result Poincar\'e conjectured,.
in \cite{B}. This was perhaps the first of Birkhoff's
many great contributions to bring him world recognition in the mathematical
community.

\begin{thm}[G. D. Birkhoff]
If $f : \A \to \A$ is an area preserving annulus homeomorphism
which preserves orientation and boundary components.  Let
$F : \tiA \to \tiA$  be a lift and let 
$F_0(x) = F(x,0),\ F_1(x) = F(x,1).$
If $p/q \in \Q$ and 
\[
\tau( F_0) \le \frac{p}{q} \le \tau( F_1)
\]
then $F$ has at least two points with translation number $p/q$
which are lifts of periodic points of $f$ with distinct 
orbits of period $q$ and rotation number $p/q + \Z$.
\end{thm}

In addition to the application which Poincar\'e had in mind this result
can be applied to billiards on a table with smooth convex boundary
as described above.  In fact, if $f : \A \to \A$ is the billiard map,
there is a lift $F : \tiA \to \tiA$ with the property that for 
all $x \in \R,\ F_0(x) = x,$ and $F_1(x) = x + 1.$  It follows that
$\tau( F_0) =0$ and $\tau( F_1) = 1$ so according to the theorem
for every $p/q \in [0,1]$ there are two periodic orbits with
translation number $p/q$ and rotation number $p/q + \Z$.  This is
true for any billiard table as long as it has a smooth strictly 
convex boundary.  (In fact, strict convexity is not really necessary;
simple convexity and smoothness of the boundary suffices.)

Figure 1 shows examples with two periodic billiard
orbits on a convex table with rotation numbers 
$\rho(x,f) = \frac{1}{5} + \Z$ and $\rho(x,f) = \frac{2}{5} + \Z.$

\begin{figure}[htb]
		\mbox{
		\epsfxsize=3.5in
		\centerline{\epsfbox{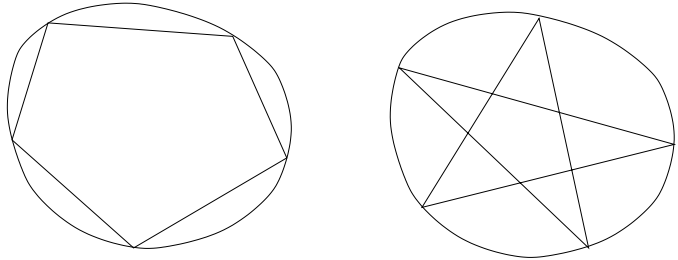}}
		}
\caption{Periodic billiard orbits with 
$\rho(x,f) = \frac{1}{5} + \Z$ and $\rho(x,f) = \frac{2}{5} + \Z$}
\end{figure}

%\begin{figure}[htb]
%\centerline{\includegraphics[width=3.5in]{b2.png}}
%\caption{Periodic billiard orbits with 
%$\rho(x,f) = \frac{1}{5} + \Z$ and $\rho(x,f) = \frac{2}{5} + \Z$}
%\end{figure}

The Poincar\'e-Birkhoff theorem has a number of generalizations.  One
which is important for our purposes asserts that we do not 
need information about the behavior of $F$ on the boundary
of $\tiA$, it suffices to have any two points $w_0, w_1 \in \tiA$
with distinct translation numbers to obtain periodic lifts
with all rational translation numbers in the interval
$[\tau( w_0, F), \tau( w_1, F)].$

\begin{thm}[Franks \cite{Fr1},\cite{Fr2}]
Suppose $f : \A \to \A$ is an area preserving annulus homeomorphism
which preserves orientation and boundary components.  Let
$F : \tiA \to \tiA$ be a lift. 
If $w_0, w_1 \in \tiA$ and  $p/q \in \Q$ satisfy
\[
\tau( w_0, F) \le \frac{p}{q} \le \tau( w_1, F)
\]
then $F$ has at least two points with translation number $p/q$
which are lifts of periodic points of $f$ with distinct 
orbits of period $q$ and rotation number $p/q + \Z$.
\end{thm}

Of course, it is natural to ask about points with irrational 
translation number.

\begin{thm}[M. Handel \cite{H1}]
Suppose that the homeomorphism $f:\A \to \A$ 
preserves orientation and boundary components and $F$ is a lift.
Then the set 
$\displaystyle \{\tau(w, F)\ |\ w \in \tiA\}$ is closed.
\end{thm}

Notice that unlike the Poincar\'e-Birkhoff Theorem, no area preserving
hypothesis is necessary for Handel's result.  However when the area
preserving hypothesis holds, we can conclude from the results above
that the set of all translation numbers of a lift $F$ forms a 
closed interval.

Finally we note that there is also a result relating 
rationality of the mean translation number to the existence of
periodic points.

\begin{thm}[Franks \cite{Fr3}]
Suppose $f : \A \to \A$ is an annulus homeomorphism
which preserves orientation and boundary components.  Let
$F : \tiA \to \tiA$ be a lift.  If $f$ preserves Lebesgue 
measure $\mu$ and 
\[
\tau_\mu( F) = \frac{p}{q}
\]
then $F$ has at least two points with translation number $p/q$
which are lifts of periodic points of $f$ with distinct 
orbits of period $q$ and rotation number $p/q + \Z$.
\end{thm}

\section{The geodesic return map on $S^2$}

In this section we discuss another application of translation and
rotation numbers,  namely the problem
of the existence of closed geodesics on the two-sphere $S^2.$

\begin{figure}[htb]
		\mbox{
		\epsfxsize=3in
		\centerline{\epsfbox{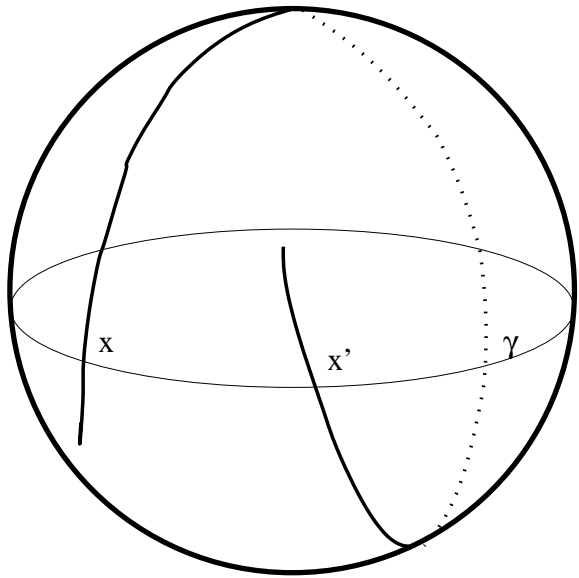}}
		}
	\caption{The Geodesic Return Map}
\end{figure}

One of the standard ways to find closed geodesics is to try to minimize
the length of a closed curve in its homotopy class.  Of course, this
cannot immediately be applied to $S^2$ since all closed curves are
homotopic to curves of arbitrarily small length.  Nevertheless 
a fascinating result of Birkhoff, using a mini-max argument, does
succeed in obtaining at least one closed geodesic for $S^2$ provided
with an arbitrary Riemannian metric.  In fact the geodesic so obtained
is {\em simple}, i.e. does not cross itself.

\begin{thm}[Birkhoff]
For any Riemannian metric on $S^2$ there is always at least one simple
closed geodesic. 
\end{thm}

In fact, Liusternik and Schnirrelman subsequently proved that there
are at least three simple closed geodesics on such a sphere.

Birkhoff realized that, like the problem of periodic billiards, the
existence of closed geodesics could often be reduced to the study
of an annulus map.

\begin{defn}
The {\em geodesic return map} \\
$F : \A \to \A$ for a simple closed geodesic $g \subset S^2$
is defined by 
\[
F(x, \theta) = (x', \theta'),
\]
where $x'$ is the second
point where $g$ is crossed by a geodesic which starts at $x$ and
forms angle $\theta$ with $g.$ $\theta'$ is the new angle.
{\rm (See Figure 2).}
\end{defn}

\begin{thm}[Birkhoff]
The geodesic return map is well-defined for a sphere if it has
positive curvature.  Whenever it is well defined it is an 
area preserving diffeomorphism.
\end{thm}

This motivates an interest in the existence of periodic orbits 
for area preserving annulus maps since there is a one-to-one 
correspondence between closed geodesics and periodic orbits of
the geodesic return map.  
The following result shows, in particular,
that whenever the geodesic return map is well-defined for a Riemannian
metric on $S^2$, there are infinitely many distinct closed geodesics
for that metric.

\begin{thm}[Franks \cite{Fr3}]\label{thm: perpoint}
Suppose $f : \A \to \A$ is an area preserving annulus diffeomorphism
which preserves orientation and boundary components.  If $f$ has
at least one periodic point, then it has infinitely many.
\end{thm}

\begin{thm}[V. Bangert \cite{B}]
When the geodesic return map is {\em not} well defined there are infinitely
many closed geodesics on $S^2$.
\end{thm}

Combining these results we have the following result which answers a 
long-standing question about the existence of closed geodesics on
$S^2.$

\begin{thm}
Every smooth Riemannian metric on $S^2$ has 
infinitely many closed geodesics.
\end{thm}

Shortly after the appearance of this result a rather different proof 
using Morse theory on the space of free loops on $S^2$ was
provided by Nancy Hingston \cite{Hi}.

Theorem~\ref{thm: perpoint} can be strengthened in several ways.
For example in \cite{FH} the following result is proved.

\begin{thm}[Franks and Handel, \cite{FH}]
Suppose 
\[
f : \A \to \A
\]
is an area preserving annulus homeomorphism
which preserves orientation and boundary components.  If $f$ has
at least one periodic point, then the periods of periodic points
are unbounded.
\end{thm}

In fact this result is valid for more general surfaces.  An area
preserving diffeomorphism  of a compact surface
is called a Hamiltonian diffeomorphism
if it is a commutator in the group of area preserving diffeomorphisms 
homotopic to the identity. In \cite{FH} it is shown that any
Hamiltonian diffeomorphism of a surface of genus $\ge 1$ must have
periodic points of arbitrarily high period.

\section{The Instability Set}\label{sec: instab}

The purpose of this section is to delve more deeply 
into the investigation of the dynamics of 
typical area preserving diffeomorphisms of annuli and spheres.

% @@@@@@@@@@@@@@@@@@@@@@@@@@@@ FIGURE @@@@@@@@@@@@@@@@@@@@@@@@@@@@@@
%\begin{figure}[htb]
%\centerline{\includegraphics[width=4in]{black.png}}
%	\caption{The Standard Map, k = 0.15}
%\end{figure}
% @@@@@@@@@@@@@@@@@@@@@@@@@@@@@@@@@@@@@@@@@@@@@@@@@@@@@@@@@@@@@@@@@@@

% @@@@@@@@@@@@@@@@@@@@@@@@@@@@ FIGURE @@@@@@@@@@@@@@@@@@@@@@@@@@@@@@
\begin{figure}[htb]
		\mbox{
		\epsfxsize=4in
		\centerline{\epsfbox{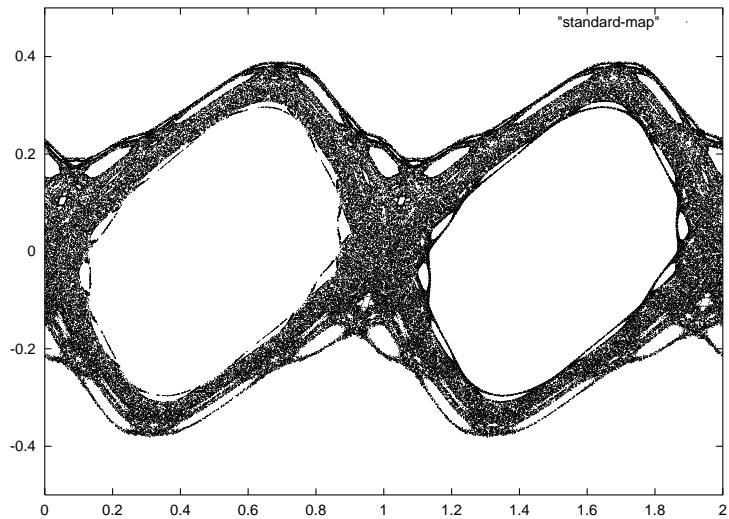}}
		}
	\caption{The Standard Map, k = 0.15}
\end{figure}
% @@@@@@@@@@@@@@@@@@@@@@@@@@@@@@@@@@@@@@@@@@@@@@@@@@@@@@@@@@@@@@@@@@@

A rather remarkable qualitative behavior is observed if one does
graphical computer studies of even simple diffeomorphisms of this
type.  The {\em standard map family} $f: \T^1 \times \R \to \T^1 \times \R$ 
is given by
\[
f(x,y) = \big (x + y + k\sin( 2\pi x), y + k\sin( 2\pi x)\big ),
\]
where $k$ is a parameter.  One checks easily that the determinant of
the Jacobian for this diffeomorphism of the open cylinder is $1$,
so it is an area preserving diffeomorphism.  If $k = 0$ this map
exhibits a rather simple dynamics, namely each circle $\T^1 \times y_0$
in the cylinder is invariant and the restriction of $f$ to this circle
is a rotation by $y_0.$

If the parameter $k$ is increased the dynamics becomes much more
complicated.  For example Figure 3 shows a computer generated image of
several thousand points on an orbit for the standard map with $k
= 0.15.$ Actually, $f(x,y)$ is periodic of period one, and what has
been drawn is the region $0 \le x \le 2,\ -0.5 \le y \le 0.5$ in order
to give a more informative picture.

The striking feature of this image is that the plotted points (a
finite set) appear to be a connected set which is riddled with holes.
The remarkable property of the image is that if another image were
created by plotting points on the orbit of another point near this set
(i.e. not in one of the holes), then the image obtained would be
visually essentially identical.  This is despite the fact that it
consists of a disjoint set of points.  As we will see, it is likely
that there is a compact connected invariant set which these plots
are approximating.

If we were to fill all the holes in this figure and take the interior
of the resulting set we would obtain an example of what Birkhoff
called a {\em ``region of instability''}.  Our aim however, is to
investigate the black set in the figure which we
will define, and refer to as an {\em instability set}.

Of particular importance in our investigation are the 
{\em stable and unstable sets} which we now define.

\begin{defn}
Suppose $f : X \to X$ is a homeomorphisms and $P \subset X$ is a
compact $f$-invariant set.  The stable set $W^s(P)$ is defined by
\[
W^s(P) := \big \{x \in X\ |\ \lim_{n \to \infty} \dist( f^n(x), P) = 0\big \},
\]
and the unstable set $W^u(P)$ is defined by
\[
W^u(P) := \big \{x \in X\ |\ \lim_{n \to -\infty} \dist( f^n(x), P) = 0\big \},
\]
where $\dist(x,P)$ is the distance from $x$ to the set $P$.
\end{defn}

If $f: S \to S$ is an area preserving diffeomorphism of  a surface
and $p \in S$ is a fixed point, then $p$ is called
{\em elliptic} if the eigenvalues of $Df_p$ are complex conjugates of
 modulus $1$ and are not equal to $\pm 1.$  The point $p$ is 
{\em hyperbolic} if the eigenvalues of $Df_p$ are real and have modulus
distinct from $1.$  Since the diffeomorphism is area preserving and
orientation preserving the product of the eigenvalues 
$ = det(Df_p) = 1.$

If $p$ is hyperbolic, then its stable set
is an injectively immersed copy of $\R,$ as is its unstable set,
and these are usually referred to as the {\em stable and unstable
manifold} of $p$.  If $p$ is a periodic point of period $n$ and
a hyperbolic fixed point of $f^n$, then $W^s(p)$ and $W^u(p)$ will
denote the stable and unstable manifold with respect to $f^n.$

We will say that an elliptic periodic point $p$ of period $n$ is {\it Moser
stable } if it admits a fundamental system of neighborhoods which are
closed disks $D$ such that $f^n|_{\partial D}$ is minimal.  A consequence
of KAM theory is that this is the typical behavior for elliptic periodic
points (see Theorem~\ref{thm: Moser generic} below).

\begin{defn}\label{moserg_def}
Suppose $f$ is an area preserving diffeomorphism
of an annular domain $\A_0$ and $\Per_q(f)$ is the set of periodic points of
$f$ with  period $q$.  Assume that for all $q$
\begin{enumerate}
\item[i)] 
the set $\Per_q(f)$ is empty or finite and all periodic points are
elliptic or hyperbolic.

\item[ii)] every elliptic periodic point in $\Per_q(f)$ is Moser stable, and
the invariant, $f$ minimal, Jordan curves bounding neighborhoods of an
elliptic periodic point $z$ have rotation numbers which are not
constant in any neighborhood of $z,$ and

\item[iii)] the intersection of any branches of hyperbolic periodic points
are transverse and any two branches of the same periodic point have
non-empty intersection.
\end{enumerate}

Then $f$ will be called {\em Moser generic.} 
\end{defn}

Combining the work of many authors one can conclude the following

\begin{thm}\label{thm: Moser generic}
If $\Sigma = S^2$ or $\A$ and $r\geq 16$, then the set of Moser generic
diffeomorphisms of $\Sigma$ is a residual subset of of the group ${\rm
Diff}^r_{\omega}(\Sigma)$ of $C^r$ area preserving diffeomorphisms.
\end{thm}

We remark that 
the high smoothness appears to be required to have all
elliptic points be Moser stable and the best published result of which
we are aware is Douady \cite{Do}, which requires $r \ge 16$. The genericity
of the properties of hyperbolic periodic points of a Moser generic
diffeomorphism follow from the Kupka-Smale theorem together with
results of Mather \cite {M3}, Pixton \cite{P}, and Robinson \cite{R} 
(cf. Theorem~\ref{thm: Pixton-Robinson} below).

The fact that the dynamics of the typical area preserving diffeomorphism
is radically different from the dynamics of the typical dissipative
diffeomorphism is amply illustrated in the following result of Mather.

\begin{thm}[Mather \cite{M3}]\label{thm: M3}
If $p$ is a hyperbolic fixed point of a Moser generic diffeomorphism of
$S^2$, then 
\[
cl( W^s(p)) = cl( W^U(p)).
\]
\end{thm}

For a Moser generic diffeomorphism there are homoclinic points
(i.e. points of $W^s(p) \cap W^U(p)$ in both components of 
$W^s(p) \setminus \{p\}$, but the reader should be cautioned that
we do not know if the homoclinic points are dense in $W^s(p).$
In any case, this common closure of the stable and unstable manifold
is a kind of organizing center for the dynamics.  Figure 3 was constructed
by plotting an orbit of the standard map, but had we picked a fixed
point and plotted points on its stable manifold we would have obtained
a visually indistinguishable picture.  All this motivates the following
definition.

\begin{defn}
Suppose $f:S^2 \to S^2$ is an area preserving diffeomorphism which is
Moser generic and $p$ is a hyperbolic point.  The 
{\em instability set} $\Lambda$ associated to $p$ is defined by
\[
\Lambda = cl( W^s(p)) = cl( W^U(p)).
\]
\end{defn}

The most fundamental open question concerning Moser generic
diffeomorphisms would be answered by a proof of the following
conjecture.

\begin{conj}\label{conj: per dense}
If $f:S^2 \to S^2$ is an area preserving diffeomorphism which is
Moser generic, then $\Per(f)$ is dense in $S^2.$
\end{conj}

Of course we would be satisfied with the proof that some other
residual subset of $\Diff^r_\omega(S^2)$ with $r$ arbitrarily large
has the property that all of its elements have a dense set of periodic
points.  This conjecture is intimately related with the question of
whether there is a $C^r$ ``closing lemma'' for area preserving $C^r$
diffeomorphisms when $r \ge 2.$ Indeed if 
Conjecture~\ref{conj: per dense} is true, then the $C^r$ ``closing lemma'' for area preserving
$C^r$ diffeomorphisms of $S^2$ holds when $r \ge 16$, because any any
$C^r$ diffeomorphism could be $C^r$ approximated by a Moser generic one
with dense periodic points.

The strongest current result in the direction of 
Conjecture~\ref{conj: per dense} the following.

\begin{thm}[Franks and LeCalvez \cite{FL}]\label{thm: W^s dense}
Suppose $f:S^2 \to S^2$ is an area preserving diffeomorphism which is
Moser generic.  Then there are countably many instability sets and their
union is dense in $S^2$.  Equivalently, if $Hyp(f)$ denotes the set of
hyperbolic periodic points, then both
\[
\bigcup_{p \in Hyp(f)} W^s(p) \text{\ \ \  and } \bigcup_{p \in Hyp(f)} W^u(p)
\]
are dense in $S^2.$
\end{thm}

As a result of this, in order to prove 
Conjecture~\ref{conj: per dense}, it would suffice to show 
that the periodic points are dense in
any instability set, or equivalently that the homoclinic points
(i.e. points in $W^s(p) \cap W^u(p)$) are dense in the instability set
associated to a hyperbolic periodic point $p.$ At present this seems
to be difficult however.  

For example, we know of no diffeomorphism 
with a hyperbolic periodic point $p$ and an open
set $U$ in which both $W^s(p)$ and $W^u(p)$ are dense but with no
periodic points and no homoclinic points in $U$.  Remarkably, however,
we don't know how to exclude this possibility even for Moser generic
diffeomorphisms! It is difficult to imagine the structure 
of any such example, however.

\section{Properties of Instability Sets}\label{sec: properties}

From its definition the instability set associated to a hyperbolic
periodic point $p$ of a Moser generic diffeomorphism $f:S^2 \to S^2$
is a compact set invariant under $f^n$ where $n$ is the period
of $p.$  It has many other striking properties
which we now summarize.

\begin{thm}[Franks and LeCalvez \cite{FL}]
Suppose $f:S^2 \to S^2$ is an area preserving diffeomorphism which is
Moser generic and $\Lambda$ is an $f$-invariant instability set.
Then
\begin{itemize}
\item
All components of $S^2 \setminus \Lambda$  are open disks.  For each such
disk $U$ there is an $n >0$ depending on $U$ such that $f^n(U) = U.$

\item
All periodic points in $\Lambda$ are hyperbolic saddles.

\item
$\Lambda$ is the closure of $W^s(q)$ for any $q \in Per(f) \cap \Lambda.$
Also $\Lambda$ is the closure of $W^u(q)$ for any $q \in Per(f) \cap \Lambda.$
\end{itemize}
\end{thm}

\begin{cor}
If $\Lambda_p$ and $\Lambda_q$ are instability sets associated with periodic 
points $p$ and $q$ respectively, then either $\Lambda_p = \Lambda_q$ or
$\Lambda_p \cap \Lambda_q = \emptyset.$
\end{cor}

We can also consider the rotation numbers of points in an instability set.
However, rotation numbers are not defined for points on $S^2$ so we must remove
two fixed points to produce an (open) annulus.  

\begin{prop}
Suppose $\Lambda$ is an instability set and $U$ is a component of $S^2 \setminus \Lambda$
satisfying $f^n(U) = U.$ Then $f^n$ has a fixed point in $U$.
\end{prop}

This is because $f^n: U \to U$ is an area preserving homeomorphism of
an open disk.  The Poincar\'e recurrence theorem asserts that almost
all points of $U$ are $f^n$ recurrent and then the Brouwer plane
translation theorem implies that $f^n$ has a fixed point.

Suppose $f:S^2 \to S^2$ is a Moser generic diffeomorphism 
and that $\Lambda$ is an $f$-invariant instability set.
If $U_0$ and $U_1$ are components of the complement of $\Lambda$ which
are invariant under $f^n$, then we may remove a fixed point of $f^n$ from
each to obtain an open $f^n$-invariant annulus in which $\Lambda$ is
essential.  Translation numbers and rotation numbers can be defined for 
an open annulus in precisely
the same way as Definition~\ref{defn: A translation number}
and Definition~\ref{defn: A rotation number}, or alternatively,
we could ``blow up'' the two fixed points rather than deleting
them and obtain a closed annulus.

We will denote by $\rho_{\Lambda}(f^n,U_0, U_1)$ the set of 
all rotation numbers of points in $\Lambda$ measured in the annulus
formed in this way.  It is not difficult to see that all these 
rotation numbers do not depend on the choice of particular fixed
points, but only on the disks $U_0$ and $U_1$ containing them.

\begin{thm}[Franks and LeCalvez \cite{FL}]
Suppose $f:S^2 \to S^2$ is a Moser generic diffeomorphism and that
$\Lambda$ is an $f$-invariant instability set. If $U_0$ and $U_1$
are any components of the complement of $\Lambda$ which are invariant
under $f^n$, then $\rho_{\Lambda}(f^n,U_0, U_1)$ is a non-trivial
closed interval.  Every rational number in $\rho_{\Lambda}(f^n,U_0,
U_1)$ is the rotation number (with respect to $U_0$ and $U_1$) of a
periodic point in $\Lambda$.
\end{thm}

In other words, for any two components of the complement of $\Lambda$
there is an entire interval of rates at which points rotate around them!

\section{Connections between frontier components}
\label{sec: connections}

An important class of area preserving annulus diffeomorphisms,
with numerous applications the class of twist maps.

\begin{defn} 
An area preserving diffeomorphism $f : \A \to \A$ is 
called a {\em twist map} if 
\[
\frac{\partial f_1}{\partial y} > 0.
\]
where $f(x,y) = (f_1(x,y), f_2(x,y)) \in \T^1 \times [0,1]$.
\end{defn} 

\begin{figure}[htb]\label{fig: twist}
		\mbox{
		\epsfxsize=4in
		\centerline{\epsfbox{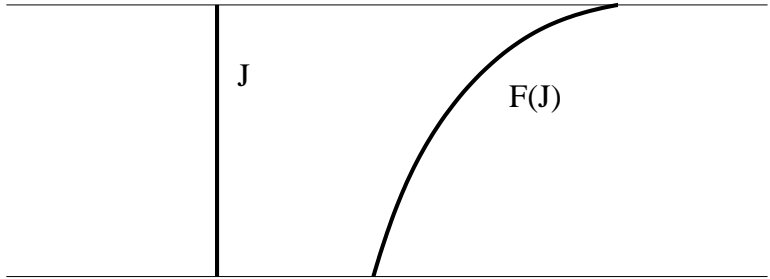}}
		}
	\caption{The lift $F$ of a Twist Map}
\end{figure}

Geometrically this means that a radial line segment
\[
J_0 = \{ (x,t)\in \A \ |\ 0 \le t \le 1\}
\]
is twisted in 
a positive sense, or
that in the lift a vertical line segment 
\[
J = \{ (x,t)\in \tiA\ |\ 0 \le t \le 1\}
\]
is bent in the positive direction (see Figure~4).

\begin{thm}[Birkhoff] \label{thm: Birkhoff}
If $f: \A \to \A$ is a twist map, any essential component of the
frontier of an invariant open set is the graph of a Lipschitz
function, in particular it is a circle.  Hence any essential component
of the frontier of an instability set is topologically a Lipschitz
circle.
\end{thm}

This is a very strong property which may be true more generally.
However, it is unknown if for a twist map or a Moser generic
diffeomorphism whether all components of the frontier must be
topological circles.  We do know that in the general (non-Moser
generic, non-twist) case the frontier of an essential open invariant
set can have components which are bad fractals.

\begin{thm}[Handel \cite{H2}]
There exists a $C^\infty$ area preserving diffeomorphism of the
annulus with an essential invariant continuum which is a
pseudo-circle.
\end{thm}

Despite the fact that we do not know if essential components of the
frontier of an instability set are topologically circles,  it was 
shown by Mather using Caratheodory's theory of prime ends that there
is always a well defined rotation number for any component 
of the frontier of an instability set of a Moser generic diffeomorphism.
Indeed Mather proves the following important result.

\begin{thm}[Mather \cite{M3}]\label{thm: Mather irrat}
Any component of the frontier of a Moser generic diffeomorphisms 
contains no periodic points and has an irrational rotation number.
\end{thm}

Mather's result actually applies somewhat more generally than the Moser
generic setting (see \cite{M3} for details).

For twist maps then, we have several interesting properties.  If
$\Lambda$ is an instability set which is essential in the annulus
There are two distinguished components, $X$ and $Y$, of the frontier
of $\Lambda$, namely those which are also the frontier of the
components of $\A \setminus \Lambda$ which contain the boundary
components of $\A.$ They are the only components of the frontier of
$\Lambda$ which are essential in $\A.$ The beautiful result of
Birkhoff, Theorem~\ref{thm: Birkhoff} above, asserts that these
are embedded Lipschitz circles.  Birkhoff also proved another striking
result about these components.  He showed that that if $U$ and $V$ are
neighborhoods of $X$ and $Y$ respectively, then there are $n, m >0$
such that $f^n(U) \cap V \ne \emptyset$ and $f^m(V) \cap U \ne
\emptyset$.  This implies that one cannot dynamically separate $X$ and
$Y$ in the sense that any open invariant set containing one of them
must contain the other in its closure.  In particular there can be no
invariant circles separating them.

This result was subsequently significantly strengthened by John Mather
who proved the following.

\begin{thm}[Mather \cite{M2}]
If $f: \A \to \A $ is a twist map and $X$ and $Y$ are the two
essential components of the frontier of $B$, then
\[
W^s(X) \cap W^u(Y) \ne \emptyset \text{ and } W^u(X) \cap W^s(Y) \ne \emptyset.
\]
\end{thm}

It would be nice to know the analogous statement for other components
of the frontier of an essential instability set.  It would also be
nice to know this with hypotheses more general than those for a twist
map.  For example, an instability set of an area preserving 
diffeomorphism of $S^2$ has no distinguished frontier components and
no distinguished co-ordinates in which to measure twist. At least in the
Moser generic case we have precisely this result.

\begin{thm}[Franks and LeCalvez \cite{FL}]\label{thm: FL}
Suppose $f: S^2 \to S^2$ is Moser generic, and $\Lambda$ is an instability
set. If $X$ and $Y$ are any two components
of the frontier of $\Lambda$, then
\[
W^s(X) \cap W^u(Y) \ne \emptyset \text{ and } W^u(X) \cap W^s(Y) \ne \emptyset.
\]
\end{thm}

In some ways the components of the frontier
of an instability set are like periodic points 
(for example, they are periodic).  Pursuing this analogy 
we should compare Theorem~\ref{thm: FL} to the following result.

\begin{thm} [ D. Pixton \cite{P}, C. Robinson \cite{R}]
\label{thm: Pixton-Robinson}
Suppose $f: S^2 \to S^2$ is Moser generic, and $\Lambda$ is an instability
set. If $p,q \in \Lambda$ are periodic points, then
\[
W^s(p) \cap W^u(q) \ne \emptyset.
\]
\end{thm}

\section{Questions about Moser generic diffeomorphisms}

The study of Moser generic diffeomorphisms (or any other residual
subset of $C^r$ area preserving diffeomorphisms with $r \ge 2$)
is a area with many more questions than answers.  The results
above as well as some additional results in \cite{FL} provide only a
start in this investigation.

Here we list some open questions about instability sets.  
They are framed in the context
of Moser generic diffeomorphisms, but answers would be very interesting
for any large class of $C^r$ diffeomorphisms with $r$ larger than $1.$
It may be useful view these questions after a further look at the
computer image of Figure 3.

For all these questions we suppose $f: S^2 \to S^2$ is a $C^r, \ r \ge
2$ Moser generic diffeomorphism preserving Lebesgue measure and let
$\Lambda$ be an instability set.  Instability sets share some
properties of hyperbolic sets (though they are never hyperbolic except
in the case of an Anosov diffeomorphism on the torus).  In particular
it may be that instability sets are topologically transitive and have dense
periodic points.

\begin{quest}[Transitivity]
Is there always a point in $\Lambda$ whose orbit is dense in $\Lambda$?
Are the periodic points of $\Lambda$ always dense in $\Lambda$?
\end{quest}

An affirmative answer to the second question 
combined with Theorem~\ref{thm: W^s dense}
would, of course, imply Conjecture~\ref{conj: per dense}.  
However, it is conceivable that this conjecture could be true without
the periodic points being dense in each instability set.

\begin{quest}[The frontier components]
Are the components of the frontier of $\Lambda$ topological
circles?  Is the restriction to such a component a minimal homeomorphism?
\end{quest}

It is quite possible that even when a component of the frontier is
a circle that the restriction of $f$ is not minimal, despite the fact
that it has irrational rotation number by the result of Mather cited
above (Theorem~\ref{thm: Mather irrat}).

\begin{quest}[The size of instability sets]
Does every instability set $\Lambda$
have positive measure?
Can $\Lambda$  have interior?  
Can the complement of $\Lambda$ have finitely many components?  
Indeed, must the complement even be non-empty?
\end{quest}

Of course, for diffeomorphisms on the torus $\T^2$ we have
the example of linear Anosov diffeomorphisms.  These are 
Moser generic and have a single instability set equal to all
of $\T^2$.  This seems likely to be a phenomenon which cannot
occur on $S^2.$

Closely related is the following question which really displays the
depth of our ignorance about Moser generic diffeomorphisms.

\begin{quest}[Elliptic Points]
If $f: S^2 \to S^2$ is a Moser generic $C^r,\ r \ge 2$ diffeomorphism 
preserving Lebesgue measure, must it have infinitely many elliptic
periodic points?  Indeed must it have at least one elliptic periodic point?
\end{quest}

\noindent
John Franks:\\
Dept. of Mathematics, Northwestern University\\
Evanston, IL 60208, USA\\
E-mail: john@math.northwestern.edu

\end{document}